\documentclass{amsart}

\usepackage{amsmath,amssymb,amsthm,mathtools}
\usepackage{graphicx}
\usepackage[hidelinks]{hyperref}
\usepackage{xcolor}

\newcommand{\qua}{\hskip 0.4em \ignorespaces}
\def\arxiv#1{\relax\ifhmode\unskip\qua\fi
  \href{https://arxiv.org/abs/#1}{\tt arXiv:\penalty-100\unskip#1}}

\def\MR#1{\relax\ifhmode\unskip\qua\fi
  \href{https://mathscinet.ams.org/mathscinet-getitem?mr=#1}{\tt MR#1}}
\def\ZB#1{\relax\ifhmode\unskip\qua\fi
  \href{https://zbmath.org/?q=an:#1}{\tt Zbl\:#1}}

\definecolor{amaranth}{rgb}{0.9, 0.17, 0.31} 
\definecolor{carrotorange}{rgb}{0.93, 0.57, 0.13} 
\definecolor{citrine}{rgb}{0.89, 0.82, 0.04} 
\definecolor{dartmouthgreen}{rgb}{0.05, 0.5, 0.06} 
\definecolor{ballblue}{rgb}{0.13, 0.67, 0.8} 
\definecolor{ceruleanblue}{rgb}{0.16, 0.32, 0.75} 
\definecolor{amethyst}{rgb}{0.6, 0.4, 0.8} 
\definecolor{amber}{rgb}{1.0, 0.75, 0.0} 
\definecolor{burlywood}{rgb}{0.87, 0.72, 0.53} 

\def\xox#1{\csname xx#1\endcsname}

\newcommand{\Z}{\mathbb Z}
\newcommand{\Q}{\mathbb Q}
\newcommand{\LambdaRing}{\Z[t, t^{-1}]}
\newcommand{\doteqeq}{\mathrel{\dot=}}

\newtheorem{theorem}{Theorem}[section]
\newtheorem{lemma}[theorem]{Lemma}
\newtheorem{proposition}[theorem]{Proposition}
\newtheorem{corollary}[theorem]{Corollary}
\newtheorem{question}[theorem]{Question}
\theoremstyle{definition}
\newtheorem{definition}[theorem]{Definition}
\theoremstyle{remark}
\newtheorem{remark}[theorem]{Remark}

\numberwithin{equation}{section}
\begin{document}

\title{Non-concordant knots with the same trace}

\author{Marc Kegel}
\address{Universidad de Sevilla,
Avda.\ Reina Mercedes s/n, 41012 Sevilla, Spain}
\email{mkegel@us.es}

\author{Kimihiko Motegi}
\address{Department of Mathematics, Nihon University,
3-25-40 Sakurajosui, Setagaya-ku, Tokyo 156--8550, Japan}
\email{motegi.kimihiko@nihon-u.ac.jp}

\date{\today}

\begin{abstract}
Miller and Piccirillo demonstrated the existence of non-concordant knots that share the same $0$--trace. 
In this short note, we generalize their result by showing that every smooth concordance class contains a knot $K$ such that, for every nonzero integer $n$, 
there exists a knot $K_n$ whose $n$--trace is orientation-preservingly diffeomorphic to that of $K$, but for which $p[K_n] \ne p[K]$ in the topological concordance group for every positive integer $p$.
Furthermore, $p[K_n]$ is distinct from $p[K']$ for any knot $K'$ obtained from $K$ by mirroring and/or reversing its orientation.
\end{abstract}


\makeatletter
\@namedef{subjclassname@2020}{\textup{2020} Mathematics Subject Classification}
\makeatother

\maketitle

\section{Introduction}

For an integer $n$, the \textit{$n$--trace} of a knot $K$ in $S^3$ is the oriented, smooth $4$--manifold $X_n(K)$
obtained by attaching an $n$--framed $2$-handle to $D^4$ along $K$.  
By definition, its boundary is the $3$--manifold $K(n)$ obtained from $S^3$ by $n$--surgery on $K$. 

There exist many different constructions of non-isotopic knots (or even infinite families of pairwise non-isotopic knots) that share the same $n$--surgery, or the same $n$--trace. We refer for example to~\cite{Brakes1980,Akbulut,Osoinach2006,AbeJongLueckeOsoinach,YasuiCorksExotic,AbeAnnulusTwists,BakerMotegi,piccirillo:shake-genus,piccirillo:Conway,ManolescuPiccirillo,TagamiRGB,qin,KP_four_surgeries,HPW}.

The following result even constructs non-concordant knots with the same $0$--trace.

\begin{theorem}[Miller--Piccirillo \cite{MillerPiccirillo}]
\label{MP}
There exist infinitely many pairs of knots $K, K'$ such that 
$X_0(K)$ and $X_0(K')$ are orientation-preservingly diffeomorphic, but $K$ and $K'$ are not smoothly concordant, even up to reversal of orientation.  \qed
\end{theorem}

The purpose of this note is to extend this result as follows. 

\begin{theorem}\label{thm:every-class}
Every smooth concordance class contains a knot $K$ such that, for each nonzero integer $n$, 
there is a knot $K_n$ whose $n$--trace is orientation preservingly diffeomorphic to the $n$--trace of $K$, 
but $p[K_n] \ne p[K]$ in the topological concordance group for all positive integers $p$.
Furthermore, $p[K_n]$ is distinct from $p[K']$ for any knot $K'$ obtained from $K$ by mirroring and/or reversing its orientation. 
\end{theorem}

\begin{remark}\label{nonzero_concordant}
The restriction $n\ne0$ in Theorem~\ref{thm:every-class} is necessary. 
If $K$ is smoothly slice and $X_0(K)$ is orientation-preservingly diffeomorphic to $X_0(K')$, 
then the trace embedding lemma implies that $K'$ is smoothly slice as well; see \cite[Corollary~1.9]{MillerPiccirillo}. Thus, the zero smooth concordance class cannot satisfy the asserted conclusion for $n=0$.
\end{remark}

Hirano \cite{Hirano2025SameNonzeroTrace} demonstrates that, for each odd integer $n$, there exists a pair of non-concordant knots with the same $n$--trace.

\medskip

Theorems~\ref{thm:every-class} and~\ref{MP} together give the following.

\begin{corollary}
    \label{cor:all-n}
For every integer $n$, there exist infinitely many pairs of knots $K, K'$ 
such that $X_n(K)$ and $X_n(K')$ are orientation-preservingly diffeomorphic, 
but $K$ and $K'$ are not smoothly concordant, even up to reversal of orientation. \qed
\end{corollary}

\begin{remark}
\label{smooth_topological}
When $n\ne0$, the pairs in Corollary~\ref{cor:all-n} are distinguished even in  topological concordance, up to reversal and concordance inverse. For $n=0$, we use only the smooth conclusion of Miller--Piccirillo. Their paper formulates topological concordance up to reversal for knots with the same oriented $0$--surgery as a conjecture \cite[Conjecture~1.15]{MillerPiccirillo}.
\end{remark}

To prove Theorem~\ref{thm:every-class}, we will prove the following more precise result. 

\begin{theorem}\label{thm:main}
There exists a ribbon knot $k$ such that, for every knot $C$ and
every nonzero integer $n$, there is a knot $K_n^C$ such that
\begin{enumerate}
\item
$X_n(K_n^C)$ is diffeomorphic to $X_n(k \# C)$, and
\item both $[K_n^C]-[k\#C]$ and $[K_n^C]+[k\#C]$ have infinite order in the  topological concordance group. The same holds with $k\#C$ replaced by its orientation reverse.
\end{enumerate}
\end{theorem}

We end this introduction with a question.

\begin{question}
For a fixed nonzero integer $n$, which pairs of smooth concordance classes admit representatives with orientation-preservingly diffeomorphic $n$--traces?
\end{question}

The answer cannot include arbitrary pairs of classes. By \cite[Theorem~1.4]{HaydenMarkPiccirillo}, a diffeomorphism $X_n(K)\cong X_n(K')$ implies $\nu(K)=\nu(K')$, except possibly when $n<0$ and $\{\nu(K),\nu(K')\}=\{0,1\}$. Since $\nu$ is a smooth concordance invariant, the classes of the unknot and the positive torus knot $T_{2,5}$, with $\nu=0$ and $2$, respectively, cannot have such representatives for any integer $n$. 

\subsection*{Conventions}
Unless explicitly stated otherwise, all manifolds, knots are smooth and oriented, 
and all diffeomorphisms are smooth and orientation-preserving. 
For two oriented manifolds $X, Y$, 
$X \cong Y$ means $X$ and $Y$ are orientation-preservingly diffeomoorphic.
Knots are considered up to isotopy.  
The symbol $-K$ denotes the mirror of the knot $K$ with its orientation reversed, so that $-K$ is the
inverse of $K$ in the concordance group. We write $K^r$ for orientation reversal without mirroring and $\overline K$ for the mirror, so that $-K=\overline K^{\,r}$. These operations all preserve the Alexander polynomial up to a unit.  
Alexander polynomials are regarded
as elements of $\LambdaRing$ up to multiplication by a unit $\pm t^n$; 
we write $f\doteqeq g$ for equality up to such a unit.

\subsection*{Acknowledgements}
We would like to thank Yume Hirano for useful discussion.
This research was conducted during KM’s visit to the IMUS (Instituto de Matemáticas de la Universidad de Sevilla) in August 2026.  
We thank the IMUS for providing office space, and KM would like to thank Marc Kegel for his warm hospitality.

\subsection{Statement on AI use}
A preliminary version of the polynomial calculations in Section 3 was developed with the assistance of ChatGPT-5.6 Sol. The authors have independently verified all the calculations by hand and take responsibility for their correctness. No presented text was created by an AI.

\subsection*{Individual grant support}
\begingroup\sloppy
MK is supported by a Ram\'on y Cajal grant (RYC2023-043251-I) and by the
project PID2024-157173NB-I00 funded by MCIN/AEI/10.13039/501100011033, ESF+
and FEDER, EU; and by a VII Plan Propio de Investigaci\'on y Transferencia
(SOL2025-36103) of the University of Sevilla.
KM has been partially supported by JSPS KAKENHI Grant Number 25K07018, 23K03110, 26K00606 and 
Nihon University research grant for 2026 Overseas Researchers. 
\endgroup

\section{The RGB construction}\label{sec:rbg}

In this section, we explain the construction of the knots in Theorems~\ref{thm:every-class} and~\ref{thm:main}.
The underlying construction is the so-called RGB link construction, 
which is a general method of creating knots that share a trace or a surgered $3$--manifold. 
There are several slightly different versions of this construction, 
where the main idea is always to describe a surgery diagram or Kirby diagram of a manifold using a $3$-component link such that two pairs of these components cancel each other out. 
See, for example, \cite{Akbulut,BakerMotegi,piccirillo:shake-genus,piccirillo:Conway,ManolescuPiccirillo,TagamiRGB,qin,KP_four_surgeries,RBG_high_dim,HPW}. Here, we will use the following variant from~\cite{DonatoneKegelLewarkTruoel}.

\begin{definition}
An oriented 3-component link $L=R\cup B\cup G$ in $S^3$ is called an \emph{RGB link} if
\begin{enumerate}
    \item $R\cup B$ is isotopic to $\mu_B\cup B$; 
    \item $R\cup G$ is isotopic to $\mu_G \cup G$, where $\mu_B$ and $\mu_G$ denote meridians of $B$ and~$G$;
    \item $R$ is labelled with a dot, and $B$ and $G$ are labelled with integers $b$ and~$g$.  
\end{enumerate}
\end{definition}

Let $L=R\cup B\cup G$ be an RGB link. Since $R$ is an unknot labelled with a dot, 
we can interpret $L$ as a Kirby diagram of the $4$-manifold $X_L$ obtained from $D^4$ by attaching a $1$-handle along $R$ (or more precisely by pushing a Seifert disk $D_R$ of $R$ into $D^4$ and removing an open tubular neighborhood of it from~$D^4$), 
and attaching $2$-handles along $B$ and $G$ with framings $b$ and~$g$, respectively. 

We write $D^4\cup h_1(R)\cup h_2(B)$ for the $4$-manifold obtained by just attaching the $1$-handle and the $2$-handle along $B$ with framing $b$ to~$D^4$. 
By assumption, there exists an isotopy such that~$B$ intersects $D_R$ transversely in a single point. 
Thus, the $1$- and the $2$-handle cancel and there exists a diffeomorphism 
\[
f_B^b\colon D^4\cup h_1(R)\cup h_2(B)\rightarrow D^4.
\]
We denote by $K_G^b$ the image of $G$ under this diffeomorphism, i.e.~$f_B^b(G)=K_G^b$. 
By interchanging the roles of $B$ and~$G$, we also obtain a knot~$K^g_B$.

If we orient the link such that $\operatorname{lk}(R,B)=\operatorname{lk}(R,G)=1$, then by construction the knots $K_B^g$ and $K_G^b$ share the same $(b+g-2\operatorname{lk}(B,G))$--trace; compare \cite[Proposition~17]{DonatoneKegelLewarkTruoel}.

\begin{figure}[ht]
  \centering
  \includegraphics[width=.55\textwidth]{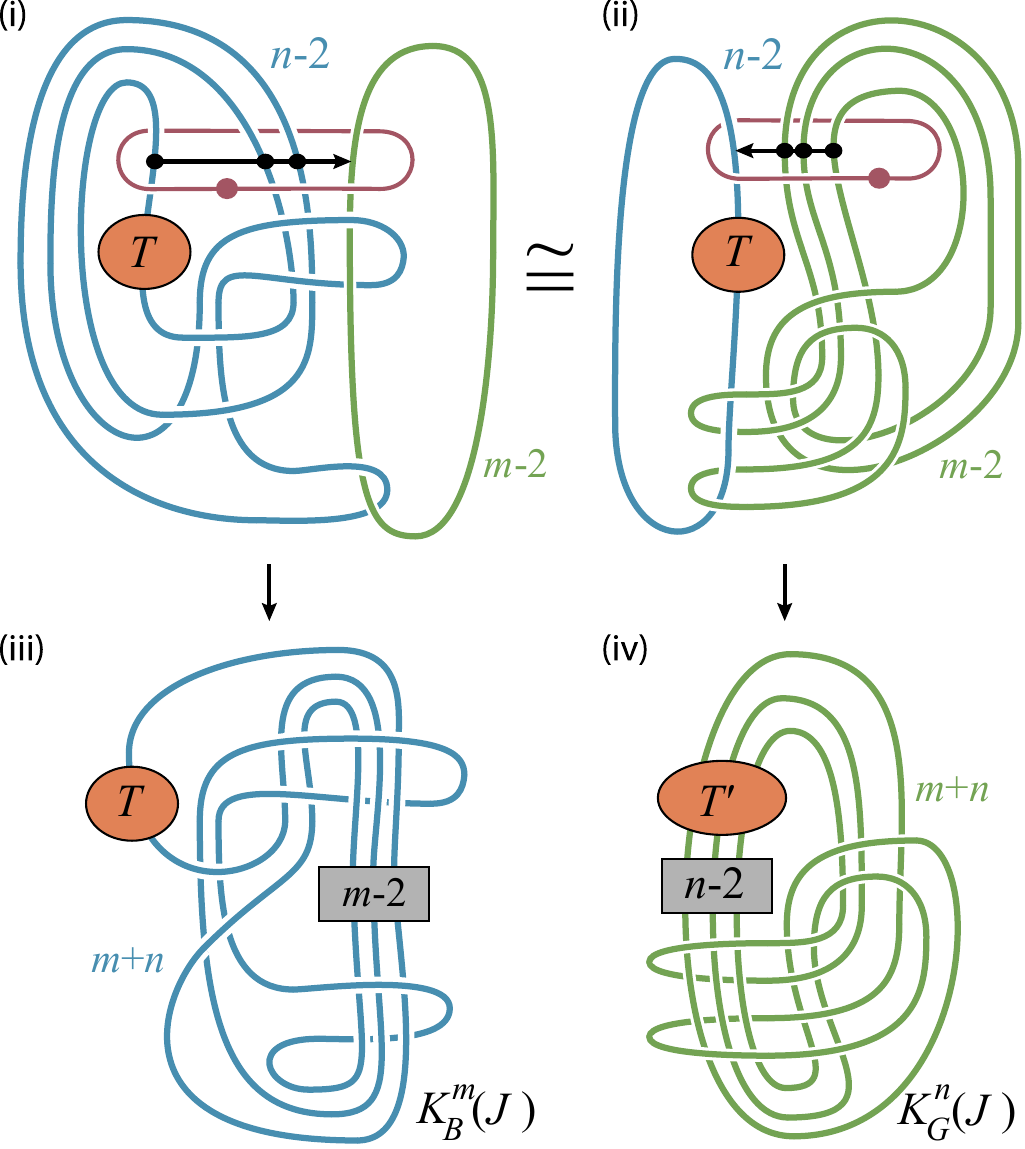}
  \caption{The RGB construction used in the proof.  The blue and green
  components in the top row have framings $n-2$ and $m-2$, respectively, and
  the red component is dotted.  The figure is reproduced from
  \cite[Figure~6]{BakerKegelMotegi} under the \textsc{cc by 4.0} license.}
  \label{fig:rbg}
\end{figure}

Consider the RGB link in Figure~\ref{fig:rbg}, taken from
\cite[Figure~6]{BakerKegelMotegi}.  
The orange disk labelled $T$ denotes a $1$--string tangle $(B^3, \tau)$ such that the closure of $\tau$ is a prescribed knot $J$. 
In addition, the orange disk labelled $T'$ represents three zero-framed parallels of the arc $\tau$. 
Here zero-framed means that the framing extends to the preferred longitude of its closure $J$: for a diagram of writhe $w$, the blackboard parallels must be corrected by $-w$ full twists. 
Write $K_B^m(J)$ and $K_G^n(J)$ for the two knots obtained after the indicated
handle slides and cancellations.  
When $J$ is the unknot, we suppress it from
the notation.

\begin{lemma}
\label{lem:companion}
For all knots $J$ and all integers $m, n$, we have
\[
  X_{m+n}\bigl(K_B^m(J)\bigr)
  \cong
  X_{m+n}\bigl(K_G^n(J)\bigr).
\]
Moreover,
\[
  K_B^m(J)=K_B^m\#J, 
\]
and for fixed $n$ there is a winding-number-one pattern $P_n$ such that
\[
  K_G^n(J)=P_n(J)
  \qquad\text{and}\qquad
  P_n(U)=K_G^n,
\]
where $U$ denotes the unknot.
\end{lemma}

\begin{proof}
Regard the red component in Figure~\ref{fig:rbg} as a dotted $1$-handle.  
The isotopy between the two diagrams in the top row, followed by the indicated handle slides, 
allows one to cancel the red component first with the green component and then, in the other presentation, with the blue component.  
With the orientations fixed above, $\operatorname{lk}(B,G)=-2$ in this diagram. Hence the remaining $2$--handle has framing $(n-2)+(m-2)-2(-2)=m+n$ in both cases.
This proves the trace diffeomorphism; it is the argument of
\cite[Section~3 and Figure~6]{BakerKegelMotegi}.

In the blue presentation, 
since the closure of the string $\tau$ is the prescribed knot $J$, we get
$K_B^m(J)=K_B^m\#J$.

In the green presentation, recall that the disk labelled $T'$ represents three zero-framed parallels of the arc $\tau$.  
The three strands have algebraic intersection number one with a meridian disk of the pattern solid torus (two strands have one orientation and the third the opposite orientation, after choosing the orientation of the core). Replacing $\tau$ by a trivial arc defines a fixed winding-number-one pattern $P_n$ with $P_n(U)=K_G^n$, and restoring $\tau$ gives $K_G^n(J)=P_n(J)$. Compare \cite[Theorem~18, Remark~19 and Figure~2]{DonatoneKegelLewarkTruoel}.
\end{proof}

For later convenience, we record the following connected-sum and satellite formulas for the Alexander polynomial of composite knots and satellite knots; see \cite[p.~118]{BurdeZieschangKnots}.

\begin{proposition}
\label{prop:alex-satellite}
In the setting of Lemma~\ref{lem:companion},
\begin{align*}
  \Delta_{K_B^m(J)}(t)
    &\doteqeq \Delta_{K_B^m}(t)\Delta_J(t),\\
  \Delta_{K_G^n(J)}(t)
    &\doteqeq \Delta_{K_G^n}(t)\Delta_J(t).
\end{align*}
\end{proposition}

\section{The Alexander-polynomial obstruction}\label{sec:polynomials}

We will distinguish the concordance classes of the knots via the Fox--Milnor condition of the Alexander polynomial. In this section, we set the necessary algebraic background of Alexander polynomials. 

For the RGB link of Figure~1 with trivial $1$--string tangle $T$,
a direct Alexander polynomial calculation gives\footnote{Up to
mirroring and reversal of orientation, the knot $K_B^2$ is
$13n_{2066}$.}
\[
\Delta_{K_B^2}(t)\doteq
f(t):=t^6-t^5-3t^4+5t^3-3t^2-t+1.
\]
The calculation immediately following Equation~$(\star\star)$
in \cite{BakerMotegi}, applied to the families in Figure~7
of that paper, gives
\[
\Delta_{K_B^m}(t)\doteq F_m(t),
\qquad
\Delta_{K_G^n}(t)\doteq F_n(t),
\]
where
\begin{equation}
F_k(t)
=(t^{-1}-2+t)(t^{-k}+t^k)
+t^{-2}-4t^{-1}+5-4t+t^2.
\end{equation}

We begin by observing some properties of the polynomial $f(t)$.

\begin{lemma}
\label{lem:q-irreducible}
The polynomial $f(t) \in \Z[t]$ is reciprocal and irreducible in $\Q[t]$.  
\end{lemma}

\begin{proof}
Reciprocity follows directly from the description of $f(t)$.  
Let
\[
\bar{f}(t)=t^6-t^5-t^3-t+1
\in\mathbb F_3[t]
\]
be the reduction of $f(t)$ modulo $3$. 
We show that $\bar{f}(t)$ is irreducible over $\mathbb F_3$. 
We use the well-known fact that
\[
t^{q^n}-t
=\prod_{\substack{g(t)\in\mathbb F_q[t]\text{ monic irreducible}\
\deg g(t)\mid n}}g(t);
\]
see, for example, Artin \cite[(6.4) Theorem~(e)]{ArtinAlgebra}. 
Thus, the irreducible factors of $t^{3^n}-t$ are precisely the monic irreducible polynomials over $\mathbb F_3$ whose degrees divide $n$.

A direct computation using the Euclidean algorithm gives
\[
\gcd(\bar{f}(t),t^3-t)=1,\qquad
\gcd(\bar{f}(t),t^9-t)=1,\qquad
\gcd(\bar{f}(t),t^{27}-t)=1.
\]
Indeed, the Euclidean algorithm gives the following remainder chains in $\mathbb F_3[t]$:
\[
\begin{aligned}
\bar f(t)&\equiv t^2+1 &&\pmod{t^3-t},\\
t^3-t&\equiv t &&\pmod{t^2+1},\\
t^2+1&\equiv1 &&\pmod t; \textrm{ and }\\
t^9-t&\equiv-t^4-t^3+1 &&\pmod{\bar f(t)},\\
\bar f(t)&\equiv t^2 &&\pmod{-t^4-t^3+1},\\
-t^4-t^3+1&\equiv1 &&\pmod{t^2}; \textrm{ and }\\
t^{27}-t&\equiv-t^5+t^4+t^2-t+1 &&\pmod{\bar f(t)},\\
\bar f(t)&\equiv1-t^2 &&\pmod{-t^5+t^4+t^2-t+1},\\
-t^5+t^4+t^2-t+1&\equiv t &&\pmod{1-t^2},\\
1-t^2&\equiv1 &&\pmod t.
\end{aligned}
\]

It follows that $\bar{f}(t)$ has no irreducible factor of degree $1, 2$, or $3$. 
Since $\bar{f}(t)$ has degree $6$, 
it must be irreducible over $\mathbb F_3[t]$.

Since the coefficient of the leading term of $f(t) \in \mathbb{Z}[t]$ is not a multiple of $3$, 
the irreducibility of $\bar{f}(t) \in \mathbb{F}_3[t]$ implies the irreducibility of $f(t)$ in 
$\Z[t]$.
Since $f(t)$ is primitive, i.e. the greatest common divisor of coefficients of $f(t)$ is $1$, 
Gauss's lemma implies that $f(t)$ is irreducible in $\Q[t]$.
\end{proof}

\begin{remark}
We may also check the irreducibility of $f(t)$ in $\mathbb{F}_3[t]$ by checking 
that $f(t)$ is not divided by monic polynomials $t^3 + a_2t^2 + a_1 t + a_0$, $t^2 + a_1 t + a_0$ and $t + a_0$, where $a_i \in \{ -1, 0, 1 \}$.
\end{remark}

\begin{lemma}\label{lem:divisibility}
$f(t)$ divides $\Delta_{K_G^n}(t)$ in
$\mathbb{Z}[t, t^{-1}]$ if and only if $n\in\{-2,2\}$.
\end{lemma}

\begin{proof}
We define $x=t+t^{-1}$ and $S_k(x)=t^k+t^{-k}$ for $k\geq 0$. 
Then $S_0(x)=2$, and $S_k(x)$ is a polynomial in $x$ of degree $k$.
Then we may write
\begin{align*}
  F_k(t) = (x-2)S_{|k|}(x)+ x^2 -4x +3 =(x-2)S_{|k|}(x)+(x-1)(x-3),\\
t^{-3}f(t)
= (t+t^{-1})^3 - (t+t^{-1})^2 - 6(t+t^{-1}) + 7 
= x^3-x^2-6x+7 
=: g(x).
\end{align*}
Since $S_2(x)=x^2-2$, these two equations give
\[
  F_2(t)=F_{-2}(t)=t^{-3}f(t).
\]

Thus $f(t)$ divides $F_{\pm2}(t)$ in $\Z[t, t^{-1}]$. 
Since $\Delta_{K_G^n}(t)\doteq F_n(t)$, it follows that $f(t)$
divides both $\Delta_{K_G^2}(t)$ and
$\Delta_{K_G^{-2}}(t)$ in $\Z[t, t^{-1}]$.

Next, we consider the converse. 
Since $g(2)=-1$ and $g(3)=7$,  $g(x) = 0$ has a real root
$\alpha\in(2,3)$.  
Choose $\lambda>1$ with
$\alpha=\lambda+\lambda^{-1}$.  
Then $f(\lambda)=0$.  
If $f(t)$ divides $F_k(t)$,
we have $F_k(\lambda)= F_2(\lambda) =0$. 
Hence, $(\alpha -2)S_{|k|}(\alpha) + (\alpha-1)(\alpha-3) = 0$ and 
$(\alpha -2)S_{2}(\alpha) + (\alpha-1)(\alpha-3) = 0$, 
which implies 
\[
  (\alpha-2)\bigl(S_{|k|}(\alpha)-S_2(\alpha)\bigr)=0.
\]
Because $\alpha>2$ and
$S_j(\alpha)=\lambda^j+\lambda^{-j}$ 
is strictly increasing for $j\geq0$, it follows that $|k|=2$.  
Applying this with $k=n$ and using
$\Delta_{K_G^n}(t)\doteq F_n(t)$ gives $n=\pm2$.
\end{proof}

\medskip

We conclude the algebraic preliminaries by recording an obstruction that follows immediately from the Fox--Milnor condition \cite{FoxMilnor}. Assume that $\varphi(t)$ is reciprocal, i.e. $\varphi(t) \doteqeq  \varphi(t^{-1})$.
If $\psi(t)\in\mathbb{Z}[t, t^{-1}]$ is divisible by $\varphi(t)$, 
then so is $\psi(t^{-1})$. 
Indeed, if $\psi(t)=\varphi(t)\eta(t)$, then
$\psi(t^{-1})=\varphi(t^{-1})\eta(t^{-1})
\doteq \varphi(t)\eta(t^{-1})$. 
We say that $\varphi(t)$ is a \textit{factor} of $\psi(t)$ of {\em exponent} $e$ if
$\psi(t)=\varphi(t)^e\eta(t)$
for some $\eta(t) \in \mathbb{Z}[t, t^{-1}]$ that is not divisible by $\varphi(t)$.

\begin{lemma}
\label{lem:fox-milnor}
Let $\phi(t)\in\mathbb{Z}[t]$ be irreducible and reciprocal,
with $\phi(0)\ne0$. 
If knots $K_0$ and $K_1$ are
 topologically concordant, then the exponent of $\varphi(t)$ in
$\Delta_{K_0}(t)\Delta_{K_1}(t)$ 
is even.
\end{lemma}

\begin{proof}
The knot $K_0\#(-K_1)$ is  topologically slice.  
The Fox--Milnor condition~\cite{FoxMilnor} says that
\[
  \Delta_{K_0}(t)\Delta_{K_1}(t)
  \doteqeq h(t)h(t^{-1})
\]
for some $h(t) \in \mathbb{Z}[t, t^{-1}]$.  
Since $\varphi(t)$ is assumed to be irreducible and  
$\varphi(t^{-1}) = \pm t^n \varphi(t)$ for some integer $n$ 
and $\mathbb{Z}[t, t^{-1}]$ is a unique factorization domain, 
the exponent of $\varphi(t)$ on the right-hand side is twice its exponent in $h(t)$.
\end{proof}

\medskip

Following Levine~\cite{Levine1969}, a Seifert matrix of a knot determines
an element of the algebraic concordance group. This association induces a well-defined epimorphism from the knot concordance group to the algebraic concordance group. Let $D$ be a Seifert
matrix of a knot $K$.

\begin{lemma}\label{lem:infinite-order}
If $f(t)$ occurs to odd exponent in $\Delta_D(t)$, 
then $D$ has infinite order in the algebraic concordance group and thus the knot $K$ has infinite order in the topological concordance group.
\end{lemma}

\begin{proof}
We first show that $f(t)$ has a root on the unit circle.
Recall that
\[
t^{-3}f(t)=g(t+t^{-1}),
\qquad
g(x)=x^3-x^2-6x+7.
\]
Since $g(1)=1$ and $g(2)=-1$, there exists
$\beta\in(1,2)$ such that $g(\beta)=0$.
Set
\[
\omega=\frac{\beta+i\sqrt{4-\beta^2}}{2}.
\]
Then $|\omega|=1$, $\omega+\omega^{-1}=\beta$, and
$f(\omega)=0$.

By Lemma~\ref{lem:q-irreducible}, $f(t)$ is irreducible over $\Q$.
Thus, since $f(\omega)=0$, the polynomial $f(t)$ is the minimal polynomial
of $\omega$ over $\Q$. 
In particular, every polynomial in $\Q[t]$
which vanishes at $\omega$ is divisible by $f(t)$.
Suppose that $f(t)$ occurs with exponent $e$ in $\Delta_D(t)$, i.e.\ there exists a polynomial $h(t)$ not divided by $f(t)$ such that
\[
\Delta_D(t)=f(t)^e h(t).
\]
Since we are in characteristic zero, the irreducible polynomial $f(t)$
has no repeated roots, so $\omega$ is a simple root of $f(t) = 0$. Moreover,
$h(\omega)\neq0$. Indeed, if $h(\omega)=0$, then the minimal polynomial
$f(t)$ of $\omega$ would divide $h(t)$, contradicting the choice of $h(t)$.
It follows that $\omega$ is a root of $\Delta_D(t) = 0$ of multiplicity
exactly $e$, which is odd by assumption.

Let $\sigma_D\colon S^1\to\Z$ denote the Levine--Tristram signature
function of $D$. We denote by $j_\omega(D)$ the jump of $\sigma_D$ at $\omega$. 
Then by~\cite[Lemma~2.1]{GaroufalidisSignature}, we have
\[
\frac{|j_\omega(D)|}{2}
\equiv
\operatorname{mult}_{\omega}(\Delta_D)
\pmod 2.
\]
Since the multiplicity on the right-hand side is odd, we get $j_\omega(D)\neq0$.
Signature jumps are additive and are invariants of algebraic
concordance. Hence, for every nonzero integer $m$,
\[
j_\omega(mD)=m\,j_\omega(D)\neq0.
\]
Thus $D$ has infinite order in the algebraic concordance group and thus $K$ has infinite order in the topological concordance group.
\end{proof}

\section{Proofs of the main results}\label{sec:proof}

Now we are ready to give the proofs of the main theorems.

\begin{proof}[Proof of Theorem~\ref{thm:main}]
We set $A=K_B^2$ and consider the ribbon knot $k=A\#(-A)$.
Let $C$ be any knot and let $n\notin\{0,4\}$.
We apply Lemma~2.2 to the pair of parameters $(2,n-2)$
and to the companion knot $J=(-A)\#C$.
Then the blue knot is
\[
K_B^2(J)=A\#(-A)\#C=k\#C,
\]
and the green knot
\[
K_n^C:=K_G^{n-2}(J)=P_{n-2}(J)
\]
has the same $n$--trace as the blue knot, since
$2+(n-2)=n$.

Recall that $\Delta_A(t)\doteq f(t)$ and thus by Proposition~2.3,
\begin{align*}
  \Delta_{k \# C}(t)& \doteqeq \Delta_k(t)\Delta_C(t) \doteqeq f(t)^2\Delta_C(t),\\
  \Delta_{K_n^C}(t)& \doteqeq \Delta_{K_G^{n-2}(J)}(t) \doteqeq \Delta_{K_G^{n-2}}(t)\Delta_J(t) 
   \doteqeq \Delta_{K_G^{n-2}}(t)f(t)\Delta_C(t).
\end{align*}

Let $e$ be the exponent of $f(t)$ in $\Delta_C(t)$.
Since $n\notin\{0,4\}$, we have $n-2\notin\{-2,2\}$.
Lemma~3.3 therefore shows that $f(t)$ does not divide
$\Delta_{K_G^{n-2}}(t)$.
The exponent of the irreducible polynomial $f(t)$ in the product $\Delta_{k \# C}(t)\Delta_{K_n^C}(t)$ is therefore
\[
  (2+e)+(1+e)=3+2e,
\]
which is odd. Thus Lemma~\ref{lem:fox-milnor} implies that $K_n^C$ is not  topologically concordant to $k\#C$. 
More strongly, for each
\[
 L\in\{k\#C,(k\#C)^r,\overline{k\#C},-(k\#C)\},
\]
the Alexander polynomial of $K_n^C\#(-L)$ has the same odd exponent $3+2e$ of $f(t)$. Lemma~\ref{lem:infinite-order} shows that $[K_n^C]-[L]$ has infinite order in both algebraic and  topological concordance.

Finally, for $n=4$, we can apply the already proved case for $n=-4$ to the knot $-C$ and then argue analogously.
\end{proof}

\begin{proof}[Proof of Theorem~\ref{thm:every-class}]
Let $C$ represent the prescribed smooth concordance class and let $k$ be the fixed ribbon knot of Theorem~\ref{thm:main}. Then $K=k\#C$ is smoothly concordant to $C$. The conclusions for every nonzero integer $n$, including their infinite-order and orientation-reversal versions, follow directly from Theorem~\ref{thm:main}.
\end{proof}

\begin{remark}
The above proof shows a slightly stronger statement. Namely, every smooth concordance class contains infinitely many distinct knots, each satisfying the conclusion of Theorem~\ref{thm:every-class} for all nonzero integers $n$. Indeed, we can choose a representative $C$ and a ribbon knot $R$ with nonconstant Alexander polynomial. Then we apply Theorem~\ref{thm:main} with $C\# \Bigl(\#_r R\Bigr)$, for $r \ge 0$. The base knots $k\#C\# \Bigl(\#_r R\Bigr)$ all represent $[C]$ and satisfy the required conclusion. Their Alexander-polynomial width strictly increase with $r$, so they are mutually distinct.
\end{remark}

\medskip

\bibliographystyle{myamsalpha}
\bibliography{HKM_trace}

\end{document}